\documentclass[11pt]{article}
\usepackage{amssymb}
\usepackage{amsmath}
\usepackage{amscd}
\usepackage{amstext}
\usepackage{amsfonts}
\usepackage{euscript}

\newcommand{\C}{\mathcal}
\newcommand{\G}{\mathfrak}
\newcommand{\N}{{\mathbb N}}

\newcommand{\F}{\mathbb {F}}
\newcommand{\X}{\mathbf X}
\newcommand{\T}{\mathbf T}

\newcommand{\gga}{\alpha}

\newcommand{\ggf}{\varphi}
\newcommand{\ggD}{\Delta}
\newcommand{\ol}{\overline}
\newcommand{\ul}{\underline}
\newcommand{\coe}{\subseteq}
\newcommand{\co}{\subset}

\newcommand{\la}{\forall}
\newcommand{\De}{\underline D}
\newcommand{\Del}{\underline{\Delta}}

\newcommand{\Spec}{\mbox{\rm Spec}}
\newcommand{\height}{\mbox{\rm ht}}
\newcommand{\Sum}[1]{{\displaystyle \sum #1}}

\DeclareMathOperator{\HS}{HS} \DeclareMathOperator{\Der}{Der}

\DeclareMathOperator{\gr}{gr}

\DeclareMathOperator{\Qt}{Qt}

\DeclareMathOperator{\ord}{ord}

\newcounter{numero}[section]
\renewcommand{\thenumero}{(\thesection .\arabic{numero})}
\newenvironment{teorema}{\medskip
\refstepcounter{numero}\noindent {\sc  \thenumero\ Theorem.}\
\it}{\vspace{1ex}\par}

\newenvironment{lema}{\medskip
\refstepcounter{numero}\noindent {\sc  \thenumero\ Lemma.}\
\it}{\vspace{1ex}\par}

\newenvironment{proposicion}{\medskip
\refstepcounter{numero}\noindent {\sc  \thenumero\ Proposition.}\
\it}{\vspace{1ex}\par}

\newenvironment{ejemplo}{\medskip
\refstepcounter{numero}\noindent {\sc  \thenumero\ Example.}\
}{\vspace{1ex}\par}

\newenvironment{demostracion}{
\noindent {\sc  Proof.}\ }{\hfill Q.E.D.\vspace{1ex}\par}

\newcommand{\numero}{\refstepcounter{numero}\noindent {\sc  \thenumero\ }}

\begin{document}

\title{Coefficient fields and scalar extension in positive characteristic}

\author{M. Fern\'{a}ndez-Lebr\'{o}n, L. Narv\'{a}ez-Macarro\thanks{Both
authors are partially supported by MTM2004-07203-C02-01 and FEDER.}\\
Department of Algebra, Faculty of Mathematics
\\ University of Seville\\ c/ Tarfia s/n, 41012 Sevilla, Spain\\ E-mail: lebron@algebra.us.es,
narvaez@algebra.us.es}
\date{}
\maketitle

\thispagestyle{empty}

\begin{abstract}
Let $k$ be a perfect field of positive characteristic, $k(t)_{per}$
the perfect closure of $k(t)$ and $A=k[[X_1,\dots,X_n]]$. We show
that for any maximal ideal $\G n$ of $A'=k(t)_{per}\otimes_k A$, the
elements in $\widehat{A'_{{\G n}}}$ which are annihilated by the
``Taylor" Hasse-Schmidt derivations with respect to the $X_i$ form a
coefficient field of $\widehat{A'_{{\G n}}}$.
\end{abstract}
\medskip

 \smallskip

\medskip

\section*{Introduction} Let
$k$ be a perfect field, $k_{(\infty)}=k(t)_{per}$ the perfect
closure of $k(t)$ and $A=k[[X_1,\dots,X_n]]$.

If $k$ is of characteristic $0$, then $k_{(\infty)} = k(t)$ and
$A(t)= A\otimes_k k(t)$ is obviously noetherian. Actually, $A(t)$ is
an $n$-dimensional regular non-local ring (see Example \ref{ejem})
whose maximal ideals have the same height ($=n$). In \cite{nar_91}
the second author proved that there is a uniform way to obtain a
coefficient field in the completions $\widehat{(A(t))_{{\G n}}}$, for
all maximal ideals $\G n$ in $A(t)$. Namely, the elements in
$\widehat{(A(t))_{{\G n}}}$ which are annihilated by the partial
derivatives $\frac{\partial}{\partial X_i}$ form a coefficient field
of $\widehat{(A(t))_{{\G n}}}$.
\medskip

In this paper, we generalize the above result to the
positive characteristic case.
\smallskip

At first sight, in positive characteristic it seems natural to
consider Hasse-Schmidt derivations instead of usual derivations (see
\cite[Theorem~3.17]{hs}), but Example \ref{ejem} shows that the
question is not so clear.
\smallskip

Consequently, in the characteristic $p>0$ case we take the scalar
extension $k \to k_{(\infty)}$ instead of $k\to k(t)$, but a new
problem appears: it is not obvious that the ring
$A_{(\infty)}=A\otimes_k k_{(\infty)}$ is noetherian. We have proved
that result in \cite{cn}.
\smallskip

The main result in this paper says that, for every maximal ideal $\G
n$ in $A_{(\infty)}$, the elements in $\widehat{(A_{(\infty)})_{{\G
n}}}$ which are annihilated by the ``Taylor" Hasse-Schmidt
derivations with respect to the $X_i$ form a coefficient field of
$\widehat{(A_{(\infty)})_{{\G n}}}$.
\medskip

Let us now comment on the content of this paper.
\smallskip

In Section 1 we introduce our basic notations and recall some
results, mainly from \cite{cn}.
\smallskip

In Section 2 we prove our main result and give the (counter)Example
\ref{ejem}.
\smallskip

In the Appendix we give a complete proof of the Normalization Lemma
for power series rings over perfect fields, which is an important
ingredient in the proof of Theorem \ref{2.3g} and that we have not
found in the literature. Our proof closely follows the proof in
\cite{abhy_64}, but the latter works only for infinite perfect
fields.

\section{Preliminaries and notations}

 All rings and algebras considered in this paper are assumed to
be commutative with unit element. If $B$ is a ring, we shall
denote by $\dim(B)$ its Krull dimension and by $\Omega (B)$ the
set of its maximal ideals. We shall use the letters $K,L,k$ to
denote fields and $\F _p$ to denote the finite field of $p$
elements, for a prime number $p$. If ${\mathfrak p} \in \Spec(B)$,
we shall denote by $\height({\mathfrak p})$ the height of
${\mathfrak p}$. Remember that a ring $B$ is said to be {\em
biequidimensional} if all its saturated chains of prime ideals
have the same length.

If $B$ is an integral domain, we denote by $\Qt(B)$ its quotient
field.

If $k$ is a ring and $B$ is a $k$-algebra, the set of all derivations
(resp. of all Hasse-Schmidt derivations) of $B$ over $k$  (cf.
\cite{has37} and \cite{mat_86}, \S 27) will be denoted by $\Der_k(B)$
(resp. $\HS_k(B)$).

\medskip
Now, we recall the notations and some results of \cite{cn} which are
used in this paper.

\medskip
For any $\F _p$--algebra $B$, we denote $B^{\sharp}:=
{\displaystyle \bigcap_{e\geq 0} B^{p^e}}$.
\smallskip

 Let $k$ be a field of characteristic $p>0$ and consider the
field extension
$$k_{(\infty)}:=\bigcup_{m\geq 0} k\left(t^{\frac{1}{p^m}}\right)\supset k(t).$$ If $k$
is perfect, $k_{(\infty)}$ coincides with the perfect closure of
$k(t)$.

For each $k$-algebra $A$, we denote $A(t):=k(t)\otimes _k A$.
 For the sake of brevity, we will write $t_m =
 t^{\frac{1}{p^m}}$ and denote
$$A_{(m)}:=A(t_m):=A\otimes_k k(t_m)=
  A(t)\otimes_{k(t)}k(t_m), \quad A_{[m]}:=A[t_m],$$
 $$A_{(\infty)}:=A\otimes_k k_{(\infty)}={\displaystyle
\bigcup_{m\geq 0} A_{(m)}},\quad A_{[\infty]}:={\displaystyle
\bigcup_{m\geq 0} A[t_m]}.$$ Each $A_{(m)}$ (resp. $A_{[m]}$) is a
free  module over $A(t)$ (resp. over $A[t]$) of rank $p^m$.
\medskip

For each prime ideal $N$
 of $A_{(\infty)}$ we denote $N_{[\infty]}:=N\cap A_{[\infty]}$,
 $N_{[m]}:=N\cap A_{[m]}$ and
 $N_{(m)}:=N\cap A_{(m)}$.
 Similarly, if $P$ is a prime ideal of $A_{[\infty]}$
we denote $P_{[m]}:=P\cap A_{[m]}$. \bigskip

\numero \label{nume} We have the following properties
\cite{cn,hs,nar_91}:
\begin{itemize}
 \item[(i)]
${\displaystyle N=\bigcup_{m\geq 0}N_{(m)}}$, ${\displaystyle
N_{[\infty]}=\bigcup_{m\geq 0} N_{[m]}}$, (resp. ${\displaystyle
P=\bigcup_{m\geq 0}P_{[m]}}$).
 \item[(ii)] $N_{(n)}\cap A_{(m)}=N_{(m)}$ and $N_{[n]}\cap A_{[m]}=N_{[m]}$ for all  $n\geq
m$ (resp. $P_{[n]}\cap A_{[m]}=P_{[m]}$ for all  $n\geq m$).
\item[(iii)] The following conditions are equivalent:
\begin{enumerate}
 \item[(a)] $N$ is maximal (resp. $P$ is maximal).
 \item[(b)] $N_{(m)}$ (resp. $P_{[m]}$) is maximal for
 some $m\geq 0$.
 \item[(c)] $N_{(m)}$ (resp. $P_{[m]}$) is maximal for
 all $m\geq 0$.
\end{enumerate}
\item[(iv)]
$\height(N)=\height(N_{[\infty]})=\height(N_{(m)})=\height(N_{[m]})$
for all $m\geq 0$. Moreover, $\dim (A_{(\infty)}) =
\dim(A_{(m)})$.
\item[(v)] \cite[Proposition~(1.4) and Theorem~(1.6)]{nar_91}
Let us assume that $A$ is noetherian and that for
every maximal ideal ${\mathfrak m}$ of $A$, the residue field
$A/{\mathfrak m}$ is algebraic over $k$. Then for every $m\geq 0$ we
have $\dim(A_{(\infty)})=\dim(A_{(m)})=\dim(A(t))$. Moreover, if $A$
is biequidimensional, universally ca\-te\-na\-rian of Krull dimension
$n$, then every maximal ideal of $A_{(\infty)}$ (or of $A_{(m)}$) has
height $n$.
\item[(vi)] \cite[Proposition~2.2]{cn} If $k$ is perfect and $B=
k[[X_1,\dots,X_n]]$, then $\Qt(B)^{\sharp} = k$.
\item[(vii)] \cite[Proposition~3.4]{cn} If $k$ is perfect, $A$ is an  integral $k$-algebra,
$K=Qt(A)$ and $K^\sharp$ is algebraic over $k$, then any prime ideal
$P\in \Spec(A_{[\infty]})$ with $P\cap k[t]=0$ and $P \cap A=0 $ is
the extended ideal of some $P_{[m_0]}$, $m_0\geq 0$.
\item[(viii)] \cite[Corollary~3.10]{cn} If $k$ is perfect, $A$ is
noetherian and for every maximal ideal ${\mathfrak m}$ of $A$, the
residue field $A/{\mathfrak m}$ is algebraic over $k$, then
$A_{(\infty)}$ is also noetherian. In particular
$k[[X_1,\dots,X_n]]_{(\infty)}$ is noetherian.
\item[(ix)] \cite[Theorem~30.6]{mat_86} Let $(R,{\G m})$ be an
equicharacteristic $n$-dimensional regular local ring containing a
quasi-coefficient field $k_0$, and $D_1,\dots,D_n\in \Der_{k_0}(R)$,
$a_1,\dots,a_n\in R$ such that $D_i(a_j)=\delta_{ij}$. Then,
$\Der_{k_0}(R)$ is a free $R$-module with basis $\{D_1,\dots,D_n\}$.
\item[(x)] \cite[Theorem~3.17]{hs}
Let $(R,{\G m})$ be an equicharacteristic $n$-dimensional regular
local ring containing a quasi-coefficient field $k_0$, and let
$\De^1, \dots ,\De^n \in \HS_{k_0}(R)$ such that their degree 1
components $\{D^1_1,\dots,D^n_1\}$ form a basis of $\Der_{k_0}(R)$.
Let $\widehat{\De}^1,\dots,\widehat{\De}^n $ be the extensions of
$\De^1,\dots,\De^n$ to $\widehat{R}$. Then, the set  $$\{
a\in\widehat{R} \mid \widehat{D}_i^j(a)=0 \quad \forall j=1,\dots,n
,\ i\geq 1 \}$$ is a coefficient field of $\widehat{R}$ (the only one
containing $k_0$).
\end{itemize}
\bigskip

\numero \label{tay} {\it Taylor expansions} (cf. \cite{mou-vi}).
\medskip

Let $n\geq 1$ be an integer. We write $\X =(X_1,\dots ,X_n)$, $\T
=(T_1,\dots T_n)$, $\X +\T =(X_1+T_1,\dots,X_n+T_n)$ and, for
$\alpha\in \N^n$,  $\X^{\alpha} =X_1^{\alpha_1} \cdots
X_n^{\alpha_n}$.

Let $A$ be the formal power series ring $k[[\X ]]$ (or the polynomial
ring $k[\X ]$). For any $f(\X )=\Sum{_{\alpha \in \N^n}}
\lambda_{\alpha} \X^{\alpha} \in A$ we
 define $\Delta^{(\alpha)}(f(\X ))$ by: $\displaystyle
f(\X +\T )=\Sum{_{\alpha \in \N^n }} \Delta^{(\alpha )}(f(\X ))
\T^{\alpha}$. One has \begin{equation} \label{eq:tay} \Delta^{(\alpha
)} (f\cdot g)=\Sum{_{\beta +\sigma =\alpha}} \Delta^{(\beta)}(f)
\Delta^{(\sigma)}(g)\end{equation} and $\alpha !\Delta^{(\alpha
)}=(\frac{\partial}{\partial X_1})^{\alpha_1}\cdots
(\frac{\partial}{\partial X_n})^{\alpha_n} $. For $i\in\N$, $1\leq
j\leq n$ and $\alpha = (0,\dots ,\stackrel{\stackrel{j}{\smile}}{i}
,\dots ,0)$ we denote $\Delta_i^j = \Delta^{((0,\dots
,\stackrel{\stackrel{j}{\smile}}{i} ,\dots ,0))}$. From
(\ref{eq:tay}) we obtain
$$ \Delta_i^j (f\cdot g)=\sum_{r+s=i} \Delta_r^j(f)
\Delta_s^j(g),$$i.e. the sequences
$\underline{\Delta}^j:=(1_A,\Delta_1^j, \Delta_2^j,\dots )$, $1\leq
j\leq n$,  are Hasse-Schmidt derivations of $A$ (over $k$) (cf.
\cite{mat_86}, \S 27).
\bigskip

Now, let us recall the following basic well known result (cf.
\cite[Propositions~5.5.3 and 5.5.6]{gro65}).
\medskip

\begin{proposicion}\label{1.1}
Let $B$ be a noetherian ring, $P$ be a prime ideal of $B[t]$ and
${\G p}=P\cap B$. Then, one of the following conditions holds:
\begin{itemize}
\item[(a)] $P={\G p}[t],\ \height (P)=\height ({\G p})$ and
$B[t]/P\simeq(A/{\G p})[t]$.
\item[(b)] $P\supset {\G p}[t],\ \height (P)=\height ({\G p})+1$
and $B[t]/P$ is an algebraic extension of $B/{\G p}$ (generated by $t
\mod  P$).
\end{itemize}
\end{proposicion}

\section{Coefficients fields and the extension $k\to k_{(\infty)}$.}

Let $k$ be a perfect field of characteristic $p>0$ and $A$ a
$k$-algebra. For every   Hasse-Schmidt derivation $\underline{\G
D}\in \HS_k(A)$, we also denote by $\underline{{\G D}} \in
\HS_{k_{(\infty)}}(A_{(\infty)})$ the extended Hasse-Schmidt
derivation. If ${\G n}\subset A_{(\infty)}$ is a maximal ideal, we
denote by $\underline{{\G D}}_{{\G n}}$ and $\widehat{\underline{{\G
D}}_{{\G n}}}$ the extended Hasse-Schmidt derivations to
$(A_{(\infty)})_{{\G n}}$ and $\widehat{(A_{(\infty)})_{{\G n}}}$,
respectively. \medskip

The following theorem generalizes Theorem 2.3 of \cite{nar_91} to the
positive characteristic case.

\begin{teorema}\label{2.3g}
Let $k$ be a perfect field of positive characteristic $p>0$,
$A=k[[X_1,\dots,X_n]]$ the power series ring and let us consider
the Hasse-Schmidt derivations $\Del^j\in \HS_k(A)$, $j=1,\dots,n$,
 defined in \ref{tay}. Then, for each maximal ideal  ${\G
n}\subset A_{(\infty)}$ the set
 $$K_0 = \left\{ a\in \widehat{(A_{(\infty)})_{\G n} }\mid
 \widehat{(\ggD_i^j)_{\G n}}(a)=0 \quad \forall j=1,\dots ,n; \ \forall
i\geq 1 \right\}$$ is a coefficient field of the complete local
ring $\widehat{(A_{(\infty)})_{{\G n}}}$.
\end{teorema}

\begin{demostracion}
  We proceed  in two steps, as in the proof of Theorem 2.3 of \cite{nar_91}: reduction to the case $n=1$ and
  treatment of this case.
\vspace{4mm}

  {\bf Step 1: the reduction.} Let us write  $P={\G n}\cap A_{[\infty]}$,
   ${\G p}={\G n}\cap A =P\cap A=P_{(m)}\cap A$. From \ref{nume}
   (iii), (iv) we know  that the ideals ${\G n}_{(m)}$ are maximal
  and $\height({\G n}_{(m)})=\height({\G n})$ for all $m\geq 0$.
By Remark (1.8)
 of \cite{nar_91}, there are only
  two possibilities for the prime ideal ${\G p}$:
\begin{enumerate}
 \item[(i)] $\height({\G p})=n$, and then
 ${\G p}=(X_1,\dots,X_n)$ and ${\G n}={\G p}^{e}$.
  \item[(ii) ] $\height({\G p})=n-1$.
\end{enumerate}
In case (i), $k_{(\infty)}$ is a coefficient field of
$(A_{(\infty)})_{\G n}$
 as well as of its completion, and $\widehat{(\ggD_i^j)_{\G n}}(k_{(\infty)})=0$
 for every $j=1,\dots ,n$, $i\geq1$. The theorem is then a
 consequence of \ref{nume} (ix), (x).
\medskip

\noindent Let us suppose we are in case (ii). By Theorem \ref{lnor}
(Normalization Lemma) there exists a new set of variables
$X_1',\dots,X_n'$ in $A$ such that
\begin{itemize}
 \item ${\G p} \cap k[[X'_1]]=(0)$,
 \item $k[[X'_1]]\hookrightarrow A/{\G p}$ is a finite extension,
and since $A/{\G p}$ is finitely generated over $k[[X'_1]]$,
$A/{\G p}$ is a finite $k[[X'_1]]$--module,
 \item $k((X'_1))\hookrightarrow \Qt(A/{\G p})$
is a separable finite extension.
\end{itemize}

Since the Hasse-Schmidt derivations of $A$ over $k$ with respect to
the variables $X_i'$ can be expressed in terms of the $\Del^j$
(\cite[Theorem~2.8]{hs}), we can suppose $X_i'=X_i$.
\medskip

Let us write $ K=A_{(\infty)}/{\G n} =
 \Qt\left(A_{[\infty]}/P\right)$, $R=A/{\G p}$, $A'=k[[X_1]]$, ${\G n}'={\G n}\cap
 A'_{(\infty)}$,
 $P'=P\cap A'_{[\infty]}={\G n}'\cap A'_{[\infty]}$ and $K'=A'_{(\infty)}/{\G n}' =
 \Qt\left(A'_{[\infty]}/P'\right)$.

 We have
$R_{[m]}={\displaystyle \frac{A_{[m]}}{{\G p}A_{[m]}}}$,
$R_{[\infty]}={\displaystyle \frac{A_{[\infty]}}{{\G
p}A_{[\infty]}}}$, $K = {\displaystyle \bigcup_{m\geq 0}
\frac{A_{(m)}}{{\G n}_{(m)}}}$ and $K' ={\displaystyle
\bigcup_{m\geq 0} \frac{A'_{(m)}}{{\G n}'_{(m)}}}$.
\medskip

Let us consider the following commutative diagram of inclusion
      $$
      \begin{array}{rlclr}
            &             & A'[t]/P'_{[0]} &             &      \\
            &\nearrow &          &\searrow &       \\
        A' &             &          &             & A[t]/P_{[0]}. \\
            &\searrow &          &\nearrow &        \\
            &             & R=A/{\G p} &             &        \\
      \end{array}
      $$
The bottom inclusions are algebraic ($R$ is a finite $A'$-module and
$P_{[0]}\cap A = {\G p}$), hence the top ones must be so. In
particular $A'[t]/P'_{[0]}$ is algebraic over $A'$, which implies
(Proposition \ref{1.1}) that $P'_{[0]}\neq 0$, then ${\G n}'_{(0)}
\neq 0$ and ${\G n}'\neq 0$. Therefore ${\G n}'$ is maximal since
$\dim(A')=1$.
\medskip

Let us show that the inclusion $K'\subset K$ is separable algebraic.
For that, it is enough to prove that the extensions
  $$ \frac{A'_{(m)}}{{\G n}'_{(m)}} \subset \frac{A_{(m)}}{{\G
  n}_{(m)}}$$are finite and separable.

 Let us write $L'=\Qt(A')=k((X_1))$, $L=\Qt(A/{\G p})$ and consider
the following diagram of field extensions
$$
\begin{array}{ccc}
  L'=\Qt(A')      & \co  &{\displaystyle
\Qt\left( \frac{A'_{[m]}}{P'_{[m]}}\right) = \frac{A'_{(m)}}{{\G n}'_{(m)}} }  \\
 \cap                &         & \cap                            \\
  L=\Qt(R) & \co &{\displaystyle \Qt
\left( \frac{A_{[m]}}{P_{[m]}}\right)= \frac{A_{(m)}}{{\G n}_{(m)}} }. \\
\end{array}
$$ \noindent  These extensions satisfy the following properties:
\begin{itemize}
 \item[i)]  $L'\co L$ is finite and
  separable. Hence, there is a primitive element $e$, $L=L'[e]$, whose
  minimal polynomial $f(X)\in L'[X]$  satisfies $f'(X)\neq 0$.
 \item[ii)]  By Proposition \ref{1.1}, the extensions
        ${\displaystyle L\co \Qt\left( \frac{A_{[m]}}{P_{[m]}}\right) }$,
$L'\co \Qt\left({\displaystyle \frac{A'_{[m]}}{P'_{[m]}}}\right)$
        are finite and generated by the class $\ol{t}$ of $t$.
\end{itemize}

\noindent Therefore,
$$ \frac{A_{(m)}}{{\G n}_{(m)}} =  \Qt\left(
\frac{A_{[m]}}{P_{[m]}}\right) = L[\ol{t}] = L'[e][\ol{t}] =
\left(\Qt\left({\displaystyle
\frac{A'_{[m]}}{P'_{[m]}}}\right)\right)[e] = \left(
\frac{A'_{(m)}}{{\G n}'_{(m)}}\right)[e]$$and the extension
$$ \frac{A'_{(m)}}{{\G n}'_{(m)}} \subset \frac{A_{(m)}}{{\G
  n}_{(m)}}$$is finite and separable for all $m\geq 0$.
Hence, $K'\subset K$ is separable algebraic.
\medskip

Let us assume that the theorem is proved for $n=1$. Then
  $$K'_0=\left\{ a \in \widehat{(A'_{(\infty)})_{{\G n}'}} \mid
    \widehat{(\ggD_i^1)_{{\G n}'}} (a)=0 \hspace{4 mm}
    \forall i\geq 1 \right\}$$
\noindent is a coefficient field of $\widehat{(A'_{(\infty)})_{{\G
n}'}}$.

We can consider $K'_0$ as a subfield of $\widehat{(A_{(\infty)})_{{\G
n}}}$ via the inclusion $\widehat{(A'_{(\infty)})_{{\G n}'}}
\hookrightarrow \widehat{(A_{(\infty)})_{{\G n}}}$. Since $K'_0
\xrightarrow{\sim} K'$ and $K'\subset K$ is separable algebraic,
we deduce that $K'_0$ is a quasi-coefficient field of
$\widehat{(A_{(\infty)})_{{\G n}}}$.

It is clear that for all $a\in K'_0 $
$$\widehat{(\ggD_i^j)_{\G n}} (a)=0 \hspace{4 mm} \forall
j=1,\dots,n, \ \forall i\geq 1.$$In particular, the
$\widehat{(\Del^j)_{\G n}}$ are Hasse-Schmidt derivations over
$K'_0$, and by \ref{nume} (ix), the $\{
 \ggD_1^1,\dots,\ggD_1^n \}$ form a basis of
$\Der_{K'_0}(\widehat{(A_{(\infty)})_{\G n}})$.

Now, by applying \ref{nume} (x), we obtain that
 $$\left\{ a \in \widehat{(A_{(\infty)})_{{\G n}}} \mid
\widehat{(\ggD_i^j)_{{\G n}}}(a)=0 \hspace{4mm} \forall
j=1,\dots,n, \ \forall i\geq 1 \right\}$$ is a coefficient field
of $\widehat{(A_{(\infty)})_{{\G n}}}$ and the theorem is proved.
\vspace{4mm}

  {\bf Step 2: the case n=1.}  Let us write $A=k[[X]]$, $L=\Qt(A)=k((X))$ and let
${\G n}$ be a maximal ideal of $A_{(\infty)}=A\otimes_k
k_{(\infty)}$. Let us denote $P={\G n}\cap  A_{[\infty]}$. By
\ref{nume} (iv), we know that
 $$\height({\G n})=\height({\G n}_{(m)})=\height(P_{[m]})=\height(P)=1.$$
\noindent As in the first step, we focus on the case ${\G n}\cap
A=(0)$ (and then $P\cap A=(0)$). Since each $A_{[m]}=A[t_m]$ is a
unique factorization domain and each $P_{[m]}$ is a prime ideal of
$A_{[m]}$ of height 1, $P_{[m]}$ is generated by an irreducible
polynomial $F_m(t_m)\in A[t_m]$ of degree $d\geq 1$ and with some
non-constant coefficient, since $P_{[m]}\cap k[t_m]=(0)$. By
irreducibility, at least one of the coefficients of $F_m(t_m)$
must be a unit, so we may assume that it is $1$.

Let us write $K=\frac{A_{(\infty)}}{{\G n}}$ and
$K_m=\frac{A_{(m)}}{{\G n}_{(m)}}$. Since $A_{(\infty)}$ and
$A_{(m)}$ are localizations of $A_{[\infty]}$ and $A_{[m]}$
respectively, it follows that
$$K=\frac{A_{(\infty)}}{{\G n}}=\Qt \left(
\frac{A_{[\infty]}}{P} \right), \quad K_m=\frac{A_{(m)}}{{\G
n}_{(m)}}=\Qt \left( \frac{A_{[m]}}{P_{[m]}}\right).
$$
\medskip

The minimal polynomial  of $\theta_m := \left(t_m \mod
P_{[m]}\right)$ over $L$ is $F_m(t_m)$. We have $K_m = L[\theta_m]$,
$$K=\bigcup_{m\geq 0} K_m = \bigcup_{m\geq 0} L[\theta_m] =
L[\theta_0,\theta_1,\theta_2,\dots],$$ where $\theta_m =
\theta_{m+1}^p$, and the inclusion $k_{(\infty)}\hookrightarrow K$ is
a $k$-morphism which sends each $t_m$ onto $\theta_m$.
\medskip

By \ref{nume} (vi), it follows that $L^{\sharp}=k((X))^{\sharp}=k$,
and we can apply \ref{nume} (vii) to conclude that there exists
$m_0\geq 0$ such that $P$ is the extended ideal of
$P_{[m_0]}=(F_{m_0}(t_{m_0}))$. Then, $P$ (resp. ${\G n}$) is the
ideal of $A_{[\infty]}$ (resp. of $A_{(\infty)}$) generated by $\mu
=F_{m_0}(t_{m_0})$. Moreover, for every $j\geq 1$, $P_{[m_0+j]}$ is
the extended ideal of $P_{[m_0]}$ and some of the coefficients of
$\mu$ is not a $p$--th power. Hence, we can take
$$F_{m_0+j}(t_{m_0+j})=F_{m_0}(t_{m_0})=F_{m_0}(t_{m_0+j}^{p^j}),\quad j\geq 1.$$

Since $k_{(\infty)}$ is perfect, the field extension
$k_{(\infty)}\co K$ is separable and, by Cohen structure theorem,
there exists a $k_{(\infty)}$--isomorphism \begin{equation}
\label{cohen}\ggf :\widehat{(A_{(\infty)})_{\G n}}
\stackrel{\sim}{\to} K[[s]]\end{equation} which induces the
identity on residue fields and  sends the regular parameter $\mu$
of $\widehat{(A_{(\infty)})_{\G n}}$ onto $s$. One has:
$$\begin{array}{rcl}
\ggf(\mu)
&=&s\\
\ggf(t_m)&=&\theta_m\\ \ggf(X)&=&X+\xi \quad \mbox{with} \quad \xi
\in (s). \\
\end{array}$$

\noindent Let us denote by
$$\Del^X=(1,\ggD_1^X,\ggD_2^X,\dots)\in \HS_k(k[[X]])$$ the
Hasse-Schmidt derivation defined in \ref{tay} and let us assume, for
the moment, that $\ggf$ satisfies the relation
\begin{equation}\label{fi} \ggf(a(X))=a(X+\xi)\coe k[[X,\xi]]\coe
K[[\xi]]\coe K[[s]]
\end{equation}
\noindent for all $a(X)\in A=k[[X]]$.
\smallskip

\noindent Then, writing $\mu =a_d(X)t_{m_0}^d+\cdots +a_0(X)$,
$$\begin{array}{l} s={\displaystyle \ggf(\mu)=\ggf\left(\sum_{r=0}^d
a_r(X)t_{m_0}^{r}\right)=\sum_{r=0}^d
\ggf\left(a_r(X)\right)\theta_{m_0}^{r}=\sum_{r=0}^d
a_r(X+\xi)\theta_{m_0}^{r} \stackrel{\ref{tay}}{=}}\\
={\displaystyle \sum_{r=0}^d\left(\sum_{i=0}^{\infty}
\ggD_i^X(a_r(X))\xi^i\right)\theta_{m_0}^{r}=\sum_{i=0}^{\infty}
\left(\sum_{r=0}^d \ggD_i^X(a_r(X))\theta_{m_0}^{r}\right)\xi^i\in
K[[\xi]]},\\
\end{array}$$
\noindent and $\xi$ must be of order one in $s$. Hence, $\xi$ is a
new variable in $K[[s]]$ and $K[[s]]=K[[\xi]]$.
\medskip

\noindent  Let us denote by $\Del'$ the unique extension of
$\Del^X$ to $K[[s]]$ through
$$A\xrightarrow{\text{scalar ext.}} A\otimes_k k_{(\infty)}\xrightarrow{\text{local.}}
(A_{(\infty)})_{\G n}\xrightarrow{\text{compl.}}
\widehat{(A_{(\infty)})_{\G n}}\xrightarrow{\varphi\ \simeq}
K[[s]],$$ which belongs to $\HS_{k_{(\infty)}}(K[[s]])$, and let
us denote by
$$\Del^{\xi}=(1,\ggD_1^{\xi},\ggD_2^{\xi},\dots )\in
\HS_K(K[[\xi]])=\HS_K(K[[s]])$$ the Hasse-Schmidt derivation
defined in \ref{tay}, this time with respect to the variable
$\xi$.
\medskip

We will show that relation (\ref{fi}) implies that
$\Del^{\xi}=\Del'$, i.e. \begin{equation} \label{2} (\ggf\circ
\ggD^X_i)(a)=(\ggD_i^{\xi}\circ \ggf)(a)\quad \forall i\geq 0,
\forall a\in k[[X]],\end{equation} and then
$$\ggf^{-1}(K)=\ggf^{-1}\left(\left\{ c\in K[[s]] \mid \ggD_i^{\xi}(c)=0,\
\la i>0 \right\}\right)=\left\{ a\in \widehat{(A_{(\infty)})_{{\G
n}}} \mid \widehat{(\ggD_i^X)_{\G n}}(a)=0,\ \la i>0 \right\}$$ is
a coefficient field of $\widehat{(A_{(\infty)})_{{\G n}}}$ and the
step 2 would be finished.
\medskip

Let $\ggf_0: A=k[[X]]\to k[[X,\xi]]$ be the local $k$-homomorphism
defined by $\ggf_0(X) = X+\xi$. Relation (\ref{fi}) says that
$\ggf(a(X)) = \ggf_0(a(X))$ for all $a(X)\in A$.
\medskip

Let $Y$ be a new variable and consider the local $k$-homomorphisms
$\delta: k[[X]] \to k[[X,Y]]$, $\varepsilon:k[[X,\xi]] \to
k[[X,\xi,Y]]$ and $\widetilde{\ggf_0}:k[[X,Y]] \to k[[X,\xi,Y]]$
defined by:
$$ \delta(X)=X+Y,\quad \varepsilon(X)=X,\quad \varepsilon(\xi) =
\xi + Y,\quad \widetilde{\ggf_0}(Y)=Y,\quad
\widetilde{\ggf_0}(X)=X+\xi.$$ Let us also consider the local
$K$-homomorphism $\Theta: K[[\xi]] \to K[[\xi,Y]]$ defined by
$\Theta(\xi) = \xi +Y$.
 Then, the
following diagram
\begin{equation*}
\begin{CD}
k[[X] @>{\ggf_0}>> k[[X,\xi]] @>{\subset}>> K[[\xi]]\\
@V{\delta}VV  @V{\varepsilon}VV @VV{\Theta}V\\
k[[X,Y]] @>{\widetilde{\ggf_0}}>> k[[X,\xi,Y]] @>{\subset}>>
K[[\xi,Y]]
\end{CD}
\end{equation*}
is commutative and we have
\begin{eqnarray*} & \displaystyle
\sum_{i=0}^{\infty} \ggD_i^{\xi}(\ggf(a)) Y^i = \Theta(\ggf(a)) =
\varepsilon (\ggf_0(a)) = \widetilde{\ggf_0} (\delta(a)) =&\\ &
\displaystyle\widetilde{\ggf_0} \left( \sum_{i=0}^{\infty}
\ggD_i^X(a) Y^i\right) = \sum_{i=0}^{\infty} \ggf_0(\ggD_i^X(a))
Y^i = \sum_{i=0}^{\infty} \ggf(\ggD_i^X(a)) Y^i&\end{eqnarray*}
for all $a\in k[[X]]$. Therefore relation (\ref{2}) is proved and
$\Del^{\xi}=\Del'$.
\medskip

The point now is to construct a $\varphi$ in (\ref{cohen}) satisfying
(\ref{fi}). We first find $\ggf(X)=X+\xi\in K[[s]]$, and for this we
state and prove the following lemma which is a generalization of
Lemma (2.3.3) of \cite{nar_91}.

\begin{lema}\label{tfi}
There exists a unique $\xi \in K[[s]]$ such that $\xi(0)=0$ of
order $1$ satisfying $$a_d(X+\xi)\theta_{m_0}^d+\cdots
+a_0(X+\xi)=s.$$
\end{lema}

\begin{demostracion}The lemma is a consequence of the implicit function theorem. Let
$G(s,\sigma)=a_d(X+\sigma)\theta_{m_0}^d+\cdots
+a_0(X+\sigma)-s\in K[[s,\sigma]]$, with
$$G(0,0)=a_d(X)\theta_{m_0}^d+\cdots
+a_0(X)=F_{m_0}(\theta_{m_0})=0.$$ \noindent We have to check that
$$\left( \frac{\partial G}{\partial \sigma}\right)
\mid_{s=\sigma=0}=a'_d(X)\theta_{m_0}^d+\cdots +a'_0(X)\neq 0 \quad
\mbox{in} \quad K.$$ \noindent Assume the contrary: then
$a'_d(X)t_{m_0}^d+\cdots +a'_0(X)$ should be a multiple of
$F_{m_0}(t_{m_0})$ in $k((X))[t_{m_0}]$ and there would be an
$\alpha\in k((X))$ such that  $$a'_r(X)=\alpha(X)a_r(X) \quad
\mbox{for every}  \quad r=0,1,\dots ,d.$$ Since some of the
coefficients $a_r$ is $1$, we deduce that $\alpha(X)=0$ and
$a'_r(X)=0$ for every $r=0,1,\dots ,d$, and then there are $b_r(X)\in
k[[X]]$ such that $a_r(X)=b_r(X^p)$. Since $k$ is perfect we conclude
that $a_r(X) = b_r(X)^p$, contradicting the fact that some of the
coefficients of $\mu$ is not a $p$-th power.

So $\left( \frac{\partial G}{\partial \sigma}\right)
\mid_{s=\sigma=0} \neq 0$, and by the implicit function theorem,
there is a unique $\xi\in K[[s]]$ such that $\xi(0)=0$ and
$G(s,\xi)=0$. Then $\xi$ has order 1 since
$$\left( \frac{\partial \xi}{\partial s}\right) (0)=\left[ \left(
\frac{\partial G}{\partial \sigma}\right) (0,0)\right] ^{-1}\neq
0.$$
\end{demostracion}
\vspace{8mm}

Let us finish the proof of Theorem \ref{2.3g}. Let $\xi \in K[[s]]$
be as in the Lemma \ref{tfi} and let us consider the local
$k$-homomorphism $$\ggf_0: \ A=k[[X]]\to k[[X,\xi]]$$ such that
$\ggf_0(X)=X+\xi$. Let us call $\ggf:A \to K[[s]]$ the composition of
$\ggf_0$ with the inclusion $k[[X,\xi]]\subset K[[\xi]]=K[[s]]$.

 We extend $\ggf$ to
$A_{(\infty)}$ by defining  $\ggf(t_m)=\theta_m\in K_m\coe K$ and
we obtain a $k_{(\infty)}$-homomorphism $\ggf:A_{(\infty)}\to
K[[s]]$ satisfying (\ref{fi}) by construction and sending
$$\mu=F_{m_0}(t_{m_0})=a_d(X)t_{m_0}^d+\cdots +a_0(X)$$ onto the
element $$a_d(X+\xi)\theta_{m_0}^d+\cdots +a_0(X+\xi)=s.$$
Therefore, the contraction of the maximal ideal $(s)$ by $\ggf$
must be ${\G n}=(\mu)$, and so we can extend $\ggf$, first to a
local $k_{(\infty)}$-homomorphism $\ggf: (A_{(\infty)})_{{\G n}}
\to K[[s]]$, and second, by completion, to $\ggf:
\widehat{(A_{(\infty)})_{{\G n}}} \to K[[s]]$

Such a $\ggf$ induces the identity map on residue fields and sends the
regular parameter $\mu=F_{m_0}(t_{m_0})$ onto $s$. Since both
local rings are regular of dimension 1, we deduce that $\gr \ggf$
is an isomorphism, and since both rings are complete, we deduce
that $\ggf:\widehat{(A_{(\infty)})_{{\G n}}}\to K[[s]]$ is a
$k_{(\infty)}$-isomorphism satisfying (\ref{fi}) as desired.
\end{demostracion}
\medskip

The following example shows that, in order to generalize Theorem
(2.3) in \cite{nar_91} to the positive characteristic case, one has
to consider the scalar extension $k \to k_{(\infty)}$ instead of $k
\to k(t)$.
\medskip

\begin{ejemplo} \label{ejem} Let $k$ be a perfect field of
characteristic $p>0$, $A=k[[X]]$ and consider the maximal ideal
${\G n} = (X^p t-1)$ in $A(t) = A\otimes_k k(t)$. Then, there is
no coefficient field of $\widehat{(A(t))_{{\G n}}}$ on which the
$\widehat{(\ggD_i^X)_{\G n}}$, $i>0$, vanish.
\medskip

Assume the contrary, i.e. there exists a coefficient field $K_0$ of
$B:= \widehat{(A(t))_{{\G n}}}$ such that $\widehat{(\ggD_i^X)_{\G
n}}(K_0)=0$ for all $i>0$, i.e. $\widehat{(\Del^X)_{\G n}}\in
\HS_{K_0}(B)$.
\smallskip

Since $\widehat{(\ggD_1^X)_{\G n}}(X)=1$, $\widehat{(\ggD_1^X)_{\G
n}}$ would be a basis of $\Der_{K_0}(B)$ by Theorem 30.6 of
\cite{mat_86}, and by Theorem 3.17 of \cite{hs} we would have the
equality
$$K_0 = \left\{ a\in B \mid \widehat{(\ggD_i^X)_{\G n}}(a)=0,\ \la i>0 \right\}.$$ In
particular $k(t)\subset K_0$.

The residue field of $B$ is
$$K = \frac{A(t)}{\G n} = \Qt\left( \frac{A[t]}{(X^p t-1)}\right) =
k[[X]][X^{-p}] = k((X)),$$where the inclusion $k(t)\hookrightarrow K$
sends $t$ to $X^{-p}$. Let $\tau:K_0 \xrightarrow{\sim} K$ be the
$k(t)$-isomorphism induced by the inclusion $K_0 \subset B$.

By Cohen structure theorem, the inclusion $K_0 \subset B$ would be
extended to an isomorphism $\psi: K_0[[s]] \xrightarrow{\sim} B$
such that $\psi(s) = X^p t-1$ ($B$ is a one dimensional complete
local noetherian local ring with parameter $X^p t-1$) and the
diagram
\begin{equation*}
\begin{CD}
K_0[[s]] @>{\psi}>{\sim}> B = \widehat{(A(t))_{{\G n}}}\\
@V{\text{res.}}VV @VV{\text{res.}}V\\
K_0 @>{\tau}>{\sim}> K
\end{CD}
\end{equation*}
is commutative.

Since $\tau^{-1}(X)$ is congruent to $X$ mod. the maximal ideal of
$B$, we deduce that $\psi^{-1}(X)$ is congruent to $\tau^{-1}(X)$
mod. $s$, i.e. $\psi^{-1}(X) = \tau^{-1}(X) + \xi$, with $\xi \in
(s)$.

On the other hand,
$$ s = \psi^{-1} (X^p t -1) = \psi^{-1}(X^p) \psi^{-1}(t) -1 =
\psi^{-1}(X)^p t - 1 = \left(\tau^{-1}(X) + \xi\right)^p t -1=
$$
$$ \left(\tau^{-1}(X)^p + \xi^p\right)t - 1 = \left(\tau^{-1}(X^p) + \xi^p\right)t -
1 = \left(t^{-1} + \xi^p\right)t - 1 = t \xi^p \in (s^p),
$$which is a contradiction.
\end{ejemplo}

\section*{Appendix: The Normalization Lemma for power series rings over perfect fields }
\renewcommand{\thesection}{A} \setcounter{numero}{0}
\label{sec-1}

In this appendix we give a proof of the normalization lema for
power series rings over an arbitrary perfect field of positive
characteristic. Our proof is an adaptation of Abhyankar's proof
\cite{abhy_64}, 23.7 and 24.5, which uses generic linear changes
of coordinates and thus requieres the field $k$ to be infinite.
\medskip

The following lemma is straightforward.

\begin{lema}\label{not1}
Let $L$ be a field of characteristic $p>0$, and let $L\subset
K=L[\gga_1,\dots ,\gga_n]$ a field extension with $\gga_i^p\in L$
for $i=1,\dots ,n$, and $[K:L]=p^{e}$. Then, there exist
$\gga_{i_1},\dots ,\gga_{i_e}$ such that $K=L[\gga_{i_1},\dots
,\gga_{i_e}]$.
\end{lema}

A series $f(X_1,\dots ,X_n)\in k[[X_1,\dots ,X_n]]$ is said to be {\em
$X_n$-distinguished} if $f(0,\dots ,0,X_n)\neq 0$.
\medskip

 The following combinatorial lemma is classical.

\begin{lema}\label{lem0}
Let $\sigma=(\sigma_1,\dots,\sigma_{n-1})\in (\N^*)^{n-1}$ and
$L_{\sigma}:\N^n \to \N$ defined by $L_{\sigma}(\gga)=\sigma_1
\gga_1+\cdots +\sigma_{n-1}\gga_{n-1}+\gga_n$ for all $\gga\in
\N^n$. Then, for each finite subset $F \subset \N^n$, there exists
a constant
 $C\geq 1$ such that the restriction $L_{\sigma}|_{F}$ is
 inyective for all $\sigma$ with  $\sigma_1\geq \sigma_2 C,$
$\sigma_2\geq \sigma_3 C$,\dots, $\sigma_{n-2}\geq \sigma_{n-1}
C$, $\sigma_{n-1}\geq C$.
\end{lema}
\begin{demostracion} The proof is standard by a double induction
on $n$ and $\sharp F$.
\end{demostracion}

\begin{lema}\label{lem1}
Let $f(X_1,\dots ,X_n)\in k[[\X]]$ be a non-zero and non-unit
formal power series. Then for $\sigma_1\gg\sigma_2\gg\cdots \gg
\sigma_{n-1}\gg 0$, the series $f(X_1+X_n^{\sigma_1},\dots
,X_{n-1}+X_n^{\sigma_{n-1}},X_n)$ is $X_n$-distinguished.
\end{lema}

\begin{demostracion}
Let us write $f(X_1,\dots ,X_n)={\displaystyle \sum_{\gga \in
\N^{n}}} f_{\gga} {X_1,\dots ,X_n}^{\gga}$ and consider the
Newton's diagram
$${\C N}(f)=\{ \gga \in \N^{n} \mid f_{\gga} \neq 0 \} \neq \emptyset
,\quad \ul 0 \notin {\C N}(f).$$ Let $F\subset {\C N}(f)$ be the
finite set of minimal elements with respect to the usual partial
ordering in $\N^{n}$. We have
 ${\C N}(f)\subset  F +\N^n$ .

\noindent By Lemma \ref{lem0}, we  obtain that $L_{\sigma}|_{F}$ is
inyective for $\sigma_1\gg \sigma_2\gg \cdots \gg \sigma_{n-1}\gg 0$,
and then the series
$$f(0+X_n^{\sigma_1},\dots ,0+X_n^{\sigma_{n-1}},X_n)={\displaystyle
\sum_{\gga \in {\C N}(f)}f_{\gga}X_n^{L_{\sigma}(\gga)}}$$has
order $\min_{\alpha\in F} L_{\sigma}(\gga)$ and is non zero.
\end{demostracion}

\begin{proposicion}\label{lem2}
Let ${\G a} \subset A=k[[X_1,\dots ,X_n]]$ be a proper ideal with
$e=\dim\left(A/{\G a}\right)$. Then there exists a
change of coordinates of the form
$$
\left\{
\begin{array}{lll}
Y_1&=&X_1+F_1(X_2^p,\dots ,X_n^p)\\ Y_2&=&X_2+F_2(X_3^p,\dots
,X_n^p)\\ \vdots& \vdots &\vdots \\ Y_{n-1}&=&X_{n-1}+F_{n-1}(X_n^p)\\
Y_n&=&X_n\\
\end{array}
\right. $$ \noindent with $F_i\in \F_p[X_{i+1},\dots,X_n]$ for
$i=1,\dots,n-1$, such that ${\G a}\cap k[[Y_1,\dots,Y_e]]=\{ 0\}$
and the extension $k[[Y_1,\dots,Y_e]]\hookrightarrow A/{\G a}$ is
finite.
\end{proposicion}
\begin{demostracion} We proceed
by induction on $n$.\\ For $n=1$: let ${\G a}$ a proper ideal of
$A=k[[X_1]]$ of height $1$. Then ${\G a}=(X_1^m)$ and $$k\subset
\frac{k[[X_1]]}{\G a}=k[\ol{X}_1] $$ is finite of rank $m$.
\medskip

Suppose now the result is true for $n-1$, and let ${\G a}$ be a
proper ideal of $A=k[[X_1,\dots ,X_n]]$. Let us take a non-zero
and non-unit formal power series $f(X_1,\dots,X_n)\in {\G a}$.

By the change
$$ \left\{
\begin{array}{l}
 Y_j=X_j-X_n^{\sigma_j},\ j=1,\dots,n-1 \\
 Y_n=X_n,\\
\end{array}
\right. $$ with $\sigma_j=\dot{p}$, $\sigma_1\gg\sigma_2\gg\cdots \gg
\sigma_{n-1}\gg 0$, and by Lemma \ref{lem1}, we deduce that the
series
$$g(Y_1,\dots,Y_{n-1},Y_n)=f(Y_1+Y_n^{\sigma_1},\dots,Y_{n-1}+Y_n^{\sigma_{n-1}},Y_n)=f(X_1,\dots ,X_n)$$
is $Y_n$-distinguished.

By Weierstrass preparation theorem we can write
$g(Y_1,\dots,Y_{n-1},Y_n)=u\cdot H$, where $u$ is a unit and
$$H=Y_n^q+a_{q-1}(Y_1,\dots,Y_{n-1})Y_n^{q-1}+\dots+a_0(Y_1,\dots,Y_{n-1}),$$
$q=\ord_{X_n}(f(X_n^{\sigma_1},\dots,X_n^{\sigma_{n-1}},X_n))\geq
1$ and $a_i(\underline{0})=0$.
\medskip

\noindent Consequently $H\in {\G a}$ and the ring extension
 $$\frac{k[[Y_1,\dots,Y_{n-1}]]}{{\G a}^c}\subseteq
\frac{k[[Y_1,\dots,Y_{n}]]}{{\G a}}=k[[Y_1,\dots,Y_{n-1}]][Y_n]$$
is finite. The proposition follows by applying induction
hypothesis to ${\G a}^c$.
\end{demostracion}
\medskip

From now on $k$ will be a perfect field of characteristic $p>0$, ${\G
p}$ a prime ideal in $A=k[[X_1,\dots ,X_n]]$, $R=A/{\G p}$,
$L=\Qt(A)=k((X_1,\dots ,X_n))$ and $K=\Qt(R)$. Let us denote $e=\dim
R$ and $\ol{a}\in R$ the class $\mod {\G p}$ of any element $a\in A$.
\medskip

The following proposition is an adaptation of (24.1) and (24.4) of
\cite{abhy_64}, which uses  Proposition \ref{lem2} instead of (23.3)
of loc. cit.
\medskip

\begin{proposicion}\label{lem3}
Under the above hypothesis, the relations
$$K=K^p[\ol{X}_1,\dots,\ol{X}_n],\  [K:K^p]=p^{e},$$hold
 and the set $\left\{\ol{X}_1^{\sigma_1} \cdots \ol{X}_n^{\sigma_n}:\
0\leq \sigma_i<p,\ i=1,\dots,n\right\}$ is a system of generators
of the extension $K^p\subset K$. Moreover, after a permutation of
variables, we have $K=K^p[\ol{X}_1,\dots,\ol{X}_e]$ and $\left\{
\ol{X}_1^{\sigma_1} \cdots \ol{X}_e^{\sigma_e} :\ 0\leq \sigma_i
<p,\ i=1, \dots,e\right\} $ \noindent is a basis of $K$ as
$K^p$-vector space.
\end{proposicion}
\begin{demostracion}
 Since $k$ is perfect, one has
$A=A^p[X_1,\dots,X_n]$, $L=L^p[X_1,\dots,X_n]$ and
$$\left\{ X_1^{\sigma_1}\cdots X_n^{\sigma_n}:\ 0\leq \sigma_1<p,\dots,0\leq
\sigma_n<p\right\}$$ is basis of $L$ (resp. of $A$) as
$L^p$-vector space (resp. as $A^p$-module). In particular $[L:L^p]
=p^n$ and $A$ is a finite $A^p$-module.

Hence, $R=R^p[\ol{X}_1, \dots, \ol{X}_n]$,
$K=K^p[\ol{X}_1,\dots,\ol{X}_n]$ and $$\left\{\ol{X}_1^{\sigma_1}
\cdots \ol{X}_n^{\sigma_n}:\ 0\leq \sigma_i<p,\
i=1,\dots,n\right\}$$ is a system of generators of the extension
$K^p\subset K$.

By Proposition \ref{lem2} we obtain a finite ring extension
$B=k[[Y_1,\dots,Y_e]]\subset R$ and then
$L_1=\Qt(B)=k((Y_1,\dots,Y_e))\subset K$ is a finite field extension.

By using Frobenius morphism one proves that $[K:L_1]=[K^p:L_1^p]$,
and from
$$[K:L_1][L_1:L_1^p]=[K:L_1^p]=[K:K^p][K^p:L_1^p]$$ we deduce that
$[K:K^p]=[L_1:L_1^p]=p^{e}$.

Finally, by Lemma \ref{not1} we know that after a permutation of
variables $$\{ \ol{X}_1^{\sigma_1} \cdots \ol{X}_e^{\sigma_e}:\ 0\leq
\sigma_i <p,\ i=1,\dots,e\}, $$ is a basis of $K$ as $K^p$-vector
space.
\end{demostracion}

\begin{teorema}\label{lnor} {\em (Normalization Lemma for power
series ring over perfect fields in positive characteristics)}
 In the situation of Proposition \ref{lem3},
 there exists a new set of variables $Y_1,\dots,Y_n \in A=k[[X_1,\dots ,X_n]]$  such that
\begin{enumerate}
 \item[(1)] ${\G p}\cap k[[Y_1,\dots,Y_e]]=\{ 0\}$.
 \item[(2)] $B=k[[Y_1,\dots ,Y_e]]\hookrightarrow R=
A/{\G p}$ is a finite ring extension.
 \item[(3)] $L_1=\Qt(B)\hookrightarrow K=\Qt(R)$ is a separable finite
extension.
\end{enumerate}
\end{teorema}
\begin{demostracion}
In view of Proposition \ref{lem3}, after a permutation of variables
$X_i$ we get $K=K^p[\ol{X}_1,\dots,\ol{X}_e]$ and $\left\{
\ol{X}_1^{\sigma_1} \cdots \ol{X}_e^{\sigma_e} :\ 0\leq \sigma_1
<p,\dots, 0\leq \sigma_e<p \right\}$ is basis of $K$ as $K^p$-vector
space.

By Proposition \ref{lem2}, there is a new set of variables
$Y_1,\dots,Y_n$ in $k[[X_1,\dots ,X_n]]$ of the form
$$ Y_j = X_j + F_j(X_{j+1}^p,\dots,X_n^p),\quad 1\leq j\leq n-1
$$and $Y_n=X_n$, with $F_j\in \F_p[X_{j+1},\dots,X_n]$, such that
${\G p}\cap k[[Y_1,\dots,Y_e]]=\{ 0\}$ and the extension
$B=k[[Y_1,\dots,Y_e]]\hookrightarrow A/{\G p}$ is finite. Hence,
$K$ is a finite field extension of $L_1=\Qt(B)$.

Since
$$
 \ol{X}_1^{\sigma_1} \cdots\ol{X}_e^{\sigma_e}
=\left(\ol{Y}_1-F_1(\ol{X}_2^p,\dots,\ol{X}_n^p)\right)^{\sigma_1}\cdots
\left(\ol{Y}_e-F_e(\ol{X}_{e+1}^p,\dots,\ol{X}_n^p)\right)^{\sigma_e},$$
$\ol{Y}_1,\dots \ol{Y}_e\in L_1=k((Y_1,\dots,Y_e))$ and
$F_j(\ol{X}_{j+1}^p,\dots,\ol{X}_n^p)=F_j(\ol{X}_{j+1},\dots,\ol{X}_n)^p
\in K^p$, we deduce that $\ol{X}_1^{\sigma_1}
\cdots\ol{X}_e^{\sigma_e}\in K^p(L_1)$ and $K=K^p(L_1)$. Therefore
$K$ is a separable finite extension of $L_1$(cf. \cite{zs_79},
Theorem 8 on p. 69).
\end{demostracion}

\end{document}